\documentclass[a4paper,12pt]{article}
\usepackage[left=2cm,right=2cm, top=2cm,bottom=2cm,bindingoffset=0cm]{geometry}
\usepackage[T2A]{fontenc}
\usepackage{graphicx}
\usepackage{cite}
\usepackage{amssymb,amsmath,latexsym}
\usepackage{caption}
\usepackage{authblk}
\usepackage{color}
\usepackage{epsfig}

\begin{document}

\title{Hyperchaos and Multistability in Nonlinear Dynamics of Two Interacting Microbubble Contrast Agents }
\author[1]{Ivan R. Garashchuk}\author[1]{Dmitry I. Sinelshchikov}
\author[2]{ Alexey O. Kazakov}
\author[1]{ Nikolay A. Kudryashov}
\affil[1]{National Research Nuclear University MEPhI, Moscow, Russia}
\affil[2]{National Research University Higher School of Economics, Nizhny Novgorod, Russia}



\maketitle

\begin{abstract}
We study nonlinear dynamics of two coupled contrast agents that are micro-meter size gas bubbles encapsulated into a viscoelastic shell. Such bubbles are used for enhancing ultrasound visualization of blood flow and have other promising applications like targeted drug delivery and noninvasive therapy.  Here we consider a model of two such bubbles interacting via the Bjerknes force and exposed to an external ultrasound field. We demonstrate that in this five-dimensional nonlinear dynamical system various types of complex dynamics can occur, namely, we observe periodic, quasi-periodic, chaotic and hypechaotic oscillations of bubbles. We study the bifurcation scenarios leading to the onset of both chaotic and hyperchaotic oscillations. We show that chaotic attractors in the considered system can appear via either Feigenbaum's cascade of period doubling bifurcations or Afraimovich--Shilnikov scenario of torus destruction. For the onset of hyperchaotic attractor we propose a new bifurcation scenario, which is based on the appearance of a homoclinic chaotic attractor containing a saddle-focus periodic orbit with its two-dimensional unstable manifold.
Finally, we demonstrate that the bubbles' dynamics can be multistable, i.e. various combinations of co-existence of the above mentioned attractors are possible. These cases include co-existence of hyperchaotic regime with any of the other remaining types of dynamics for different parameter values. Thus, the model of two coupled gas bubbles provide a new examples of physically relevant system with multistable hyperchaos.

\end{abstract}

\section{Introduction}

In this work we study nonlinear dynamics of two coupled encapsulated gas bubbles in a liquid, which are driven by an external periodic pressure field. Investigation of oscillations of such bubbles is of interest due to their applications as contrast agents for ultrasound visualization \cite{Szabo,Goldberg,Hoff} and future possible applications for noninvasive therapy and targeted drug delivery \cite{Klibanov2006,Coussios2008}. Depending on applications, different types of bubbles' dynamics can be either beneficial or undesirable (see, e.g. \cite{Hoff,Carroll2013}). Therefore, it is important to study the variety of possible dynamical regimes and how the dynamics of the bubbles depend on the both control parameters and initial conditions.


Typically, mathematical models of a single microbubble contrast agent are one-dimensional non-autonomous oscillators based on the Rayleigh--Plesset equation and its generalizations (see, e.g. \cite{Plesset,Doinikov2011} and references therein). Later, these models were extended to the coupled Rayleigh--Plesset equations which take into account bubble-bubble interactions via the Bjerknes force (see \cite{Takahira1995,Mettin1997,Ida2002,Doinikov2004,Dzaharudin2013} and references therein). From a mathematical point of view, these models are described by a system of coupled nonlinear oscillators, with external periodic force, and, thus, various types of dynamics can be observed in them. Despite the great interest, there are only few works devoted to studying of nonlinear dynamics of gas bubbles. For example, nonlinear dynamics of a single bubble described by one of the Rayleigh--Plesset-like models was studied in \cite{Parlitz1990,Behina2009,Macdonald2006,Carroll2013,Garashchuk2018,Garashchuk2018a}, where it was shown that oscillations of a single bubble can be either regular or chaotic and routes to the corresponding attractors were studied. Some bifurcations of two and $N$ coupled bubbles were studied in \cite{Takahira1995,Macdonald2006,Dzaharudin2013}. However, in \cite{Takahira1995} unencapsulated bubbles were considered, while in \cite{Macdonald2006,Dzaharudin2013} an inappropriate model of the shell was investigated (see discussion in \cite{Doinikov2011}).

After some simplifications, the dynamics of two coupled bubbles is described by a system of five ordinary differential equations. In this work we show that, in addition to regular and chaotic regimes which are typical for models of one bubble \cite{Parlitz1990,Behina2009,Macdonald2006,Carroll2013,Garashchuk2018,Garashchuk2018a}, this system exhibits quasiperiodic and, what is more interesting, hyperchaotic  types of motion. While there are many known examples of multi-dimensional nonlinear systems demonstrating quasiperiodic and chaotic dynamics with two or more positive Lyapunov exponents, to the best of our knowledge, neither quasiperiodic nor hyperchaotic oscillations of two coupled bubbles have been studied previously.

Moreover, in this work we propose a bifurcation scenario for the onset of hyperchaotic behavior with two positive Lyapunov exponents. The key part of this scenario is the appearance of a \textit{homoclinic chaotic attractor} containing a saddle-focus periodic orbit with a two-dimensional unstable manifold, i.e. such an orbit which has a pair of complex conjugated multipliers with positive real parts while all the other multipliers have negative real part. Recall, that the chaotic attractor is called homoclinic, if it contains a saddle periodic orbit \cite{GonGonPhD16} together with its unstable invariant manifold.

Trajectories on a homoclinic attractor can pass arbitrarily close to the saddle orbit belonging to it. The dynamics near this saddle orbits and, as a result, on the whole homoclinic attractor, significantly depends on the multipliers of the corresponding saddle orbit \cite{GonGonKazTur2014}. In particular, in the small neighborhood of a homoclinic attractor containing a saddle-focus periodic orbit with two-dimensional unstable manifold,  two-dimensional areas are expanded and Lyapunov exponents on the whole attractor ``can feel'' this expansion. As a result, two Lypaunov exponents become positive. Note, that homoclinic attractors of this type were called in \cite{GonGonKazTur2014}, \cite{GonGonShil2012}, \cite{GonGonPhD16} \textit{ discrete Shilnikov attractors}\footnote{This terminology arises to the paper \cite{Shil86} where Shilnikov, for one-parametric families of three-dimensional flow systems, proposed a universal bifurcation scenario leading to the birth of spiral attractor containing saddle-focus focus equilibrium together with its two-dimensional unstable manifold. In the papers \cite{GonGonShil2012}, \cite{GonGonKazTur2014} this scenario was transferred to the case of one-parametric families of three-dimensional maps (or four-dimensional flows, if the corresponding Poincar\'e cross-section is considered).}.

Another interesting property of the considered system is its multistability. It has recently been shown \cite{Garashchuk2018,Garashchuk2018a} that the dynamics of a single encapsulated bubble can be multistable, i.e. several attractors can co-exists at the same values of parameters. Thus, it is interesting to understand whether multistability persists in the dynamics of coupled bubbles or even more new multistable states can occur.  We show that the dynamics of coupled bubbles is multistable and various attractors can co-exist, including periodic and hyperchaotic attractors. Thus, we demonstrate a possibility of hidden hyperchaos in the dynamics of two interacting gas bubbles, which is a new example of physically relevant five-dimensional dynamical system with hidden hyperchaos.



The rest of this work is organized as follows. In Sec. \ref{sec:eqsys} we present the governing system of equations for the dynamics of two coupled bubbles and discuss some of its properties. 
In Sec. \ref{sec:3} we present a two-parametric chart of the Lyapunov exponents for the considered system and discuss possible types of dynamics. 
 In Sec. \ref{sec:4} we focus on scenarios leading to the onset of both chaotic and hyperchaotic oscillations of coupled bubbles. We demonstrate that a chaotic attractor can appear via either Feigenbaum's cascade of period-doubling bifurcations or Afraimovich-Shilnikov scenario of the destruction of invariant tori. We also propose a new phenomenological scenario for the onset of hyperchaotic oscillations as well. In Sec. \ref{sec:multistability} we discuss a possibility of co-existence of several attractors in the systems under consideration and point out that chaotic attractor can coexist with hyperchaotic one in this system. In the last section we briefly discuss our results specially marking that the proposed scenario of the onset of hyperchaotic attractors should be also typical for other multi-dimensional systems demonstrating hyperchaotic behavior.

%

\section{Main system of equations} \label{sec:eqsys}

\begin{figure}[!ht]
\center{\includegraphics[width=0.6\linewidth]{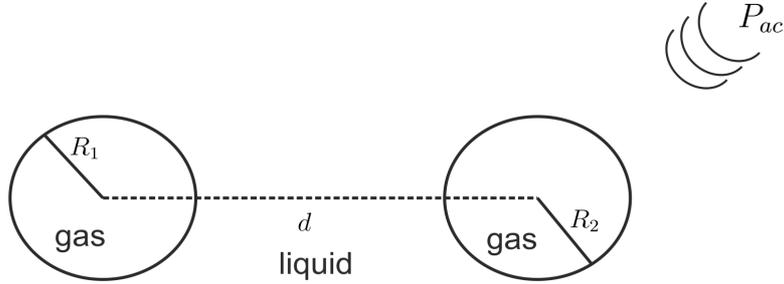} }
\caption{{\footnotesize Schematic picture of two interacting bubbles oscillating in a liquid under the influence of an external pressure field.
}}
\label{fig:A}
\end{figure}


In this section we consider a model for the description of oscillations of two interacting gas bubbles in a liquid. Essentially, this model is formed by two generalized Raleigh--Plesset equations that are coupled via the Bjerknes forces (see, e.g. \cite{Takahira1995,Mettin1997,Ida2002,Doinikov2004,Dzaharudin2013,Macdonald2006}). In this work we take into account liquid's compressibility in accordance with the Keller--Miksis model \cite{Keller1980}, liquid's viscosity on the gas-liquid interface, surface tension and the impact of bubbles' shells, which is described by the de--Jong model \cite{deJong1992,Marmottant2005}.  We also suppose that bubbles are exposed to the external periodic pressure field. Under these assumptions the governing system of equations for oscillations of two coupled bubbles is
\begin{equation}\label{eq:eq1}
  \begin{gathered}
    \left(1-\frac{\dot{R_{1}}}{c}\right)R_{1} \ddot{R_{1}} +\frac{3}{2}\left(1-\frac{\dot{R_{1}}}{3c}\right)\dot{R}_{1}^{2}=\frac{1}{\rho}\left[1+\frac{\dot{R}_{1}}{c}+\frac{R_{1}}{c}\frac{d}{dt}\right]P_{1}-\frac{d}{dt}\left(\frac{R_{2}^{2}\dot{R}_{2}}{d}\right), \vspace{0.1cm}\\ 
     \left(1-\frac{\dot{R_{2}}}{c}\right)R_{2} \ddot{R_{2}} +\frac{3}{2}\left(1-\frac{\dot{R_{2}}}{3c}\right)\dot{R}_{2}^{2}=\frac{1}{\rho}\left[1+\frac{\dot{R}_{2}}{c}+\frac{R_{2}}{c}\frac{d}{dt}\right]P_{2}
     -\frac{d}{dt}\left(\frac{R_{1}^{2}\dot{R}_{1}}{d}\right), 
  \end{gathered}
\end{equation}
where
\begin{equation*}
P_{i}=\left(P_{0}+\frac{2\sigma}{R_{i0}}\right)\left(\frac{R_{i0}}{R_{i}}\right)^{3\gamma}-\frac{4\eta_L \dot{R}_{i}}{R_{i}}-\frac{2\sigma}{R_{i}}-P_{0}-4\chi\left(\frac{1}{R_{i0}}-\frac{1}{R_{i}}\right)-4\kappa_{S}\frac{\dot{R}_{i}}{R_{i}^{2}}-P_{ac}\sin(\omega t), \quad i=1,2.
\end{equation*}

Here $R_{1}$ and $R_{2}$ are radii of bubbles, $d$ the distance between the centers of bubbles, $P_{\mathrm{stat}}$ is the static pressure, $P_v$ is the vapor pressure, $P_0 = P_{\mathrm{stat}} - P_v$, $P_{ac}$ is the magnitude of the external pressure field and $\omega$ is its cyclic frequency, $\sigma$ is the surface tension, $\rho$ is the density of the liquid, $\eta_{L}$ is the viscosity of the liquid, $c$ is the speed of sound in the liquid, $\gamma$ is the polytropic exponent, $\chi$ and $\kappa_{s}$ denote the shell elasticity and shell surface viscosity respectively.


It can be easily seen that by means of simple transformations equations \eqref{eq:eq1} can be rewritten in the form of a five-dimensional system of ordinary differential equations in terms of the following dependent variables $R_1, R_2, \dot R_1, \dot R_2$ and $\theta \in [0, 2 \pi]$. We perform all our numerical experiments exactly with this five-dimensional system. However due to its cumbersome form we do not present it here.


In what follows, we assume that $P_{ac}$, $\omega$ and $d$ are treated as the 
 control parameters and the remaining parameters are fixed as follows: $P_v = 2.33$ kPa, $\sigma = 0.0725$ N/m, $\rho= 1000$ kg/m$^3$, $\eta_{L} = 0.001$ Ns/m$^3$, $c = 1500$ m/s, $\gamma = 4/3$, $\chi = 0.22$ N/m and $\kappa_S = 2.5\cdot 10^{-9}$ kg/s. These values of the parameters correspond to the adiabatic oscillations of two interacting SonoVue contrast agents with equilibrium radii $R_{i0}=1.72 \mu$m \cite{Tu2009}.


We suppose that the equilibrium radii of bubbles are the same (i.e., $R_{10}=R_{20}=R_{0}$), because the injected ensemble of contrast agents is assumed to consist of bubbles of the same characteristics. Thus, system \eqref{eq:eq1} is symmetric with respect to  the following change of variables: $\left( R_{1}, \dot{R}_1 \leftrightarrow R_2, \dot{R}_2 \right)$. This symmetry leads to several conclusions.

First, there always exists a family of symmetrical solutions, for which $\forall t > 0 \, R_1(t) = R_2(t)$ if $R_1(t=0) = R_2(t=0)$ and $\dot{R}_1(t=0) = \dot{R}_2(t=0)$, which correspond to fully synchronous in-phase oscillations of bubbles. Some of these solutions can be asymptotically stable (attractive). Since all symmetrical solutions lie in the invariant manifold $R_1 = R_2, \dot{R}_1 = \dot{R}_2$ and system \eqref{eq:eq1} restricted to this manifold is three-dimensional and volume-contracting, such regimes can be of two possible types: periodic and simply chaotic (with only one positive Lyapunov exponent). In other words, synchronous oscillations of two bubbles can be either periodic or chaotic, for more details see the next section.

Second, asymptotically stable regimes (attractors) can exist outside the invariant manifold $R_1 = R_2, \dot{R}_1 = \dot{R}_2$. These regimes correspond to asynchronous oscillations of bubbles and, in contrast to synchronous ones, they can be of four possible types. In addition to periodic and simply chaotic regimes, asynchronous oscillations can be also quasiperiodic and even hyperchaotic. Moreover, the existence of such asymmetrical regimes lead to the trivial multistability in the system. Indeed, for each asynchronous regime passing through point $(R_1, R_2, \dot{R}_1, \dot{R}_2)$ there exist the symmetrical one passing through $(R_2, R_1, \dot{R}_2, \dot{R}_1)$. In Sec.~\ref{sec:multistability} we show that in addition to this simple multistability, system \eqref{eq:eq1} exhibits more complex types of this phenomenon. There we discuss coexistence of different attractors which cannot be obtained by simple interchanging of variables $R_1 \leftrightarrow R_2, \dot{R}_1 \leftrightarrow \dot{R}_2$  at the same values of the control parameters. Below we will refer to the term 'multistability' describing a situation of coexistence of several substantially different attractors, that are not simply symmetrical with respect to the manifold $R_1 = R_2, \dot{R}_1 = \dot{R}_2$.  Also note that coexistence of two synchronous periodic attractors both lying in the manifold $R_1 = R_2, \dot{R}_1 = \dot{R}_2$ is possible.

We perform all calculations in the following non-dimensional variables $R_{i}=R_{0}r_{i}$, $t=\omega_{0}^{-1}\tau$, where $\omega_{0}^{2}=3\kappa P_{0}/(\rho R_{10}^2)+2(3\kappa-1)\sigma/R_{10}+4\chi/R_{10}$ is the natural frequency of bubble oscillations. The non-dimensional bubbles speeds are given by $u_{i}=dr_{i}/d\tau=\dot{R}_{i}/ (R_{0}\omega_{0})$. We use the fourth-fifth order Runge--Kutta method \cite{Cash1990} for finding numerical solutions of the Cauchy problem for \eqref{eq:eq1}. For the calculations of the Lyapunov spectra we use the standard algorithm by Bennetin \cite{Benettin1980}. Poincar\'e cross sections are constructed at every period of the external pressure field.

\section{Variety of dynamical regimes in the model} \label{sec:3}

In this section we demonstrate the diversity of possible dynamical regimes in system \eqref{eq:eq1} and show that depending on the values of the control parameters $d/R_0$ and $P_{ac}$ two bubbles can exhibit periodic, quasiperiodic, chaotic or hyperchaotic oscillations. Moreover, due to the symmetry of the system, periodic and chaotic regimes can be either synchronous or asynchronous. Hyperchaotic and quasiperiodic oscillations cannot be synchronous.

\begin{figure}[!ht]
\center{\includegraphics[width=1.0\linewidth]{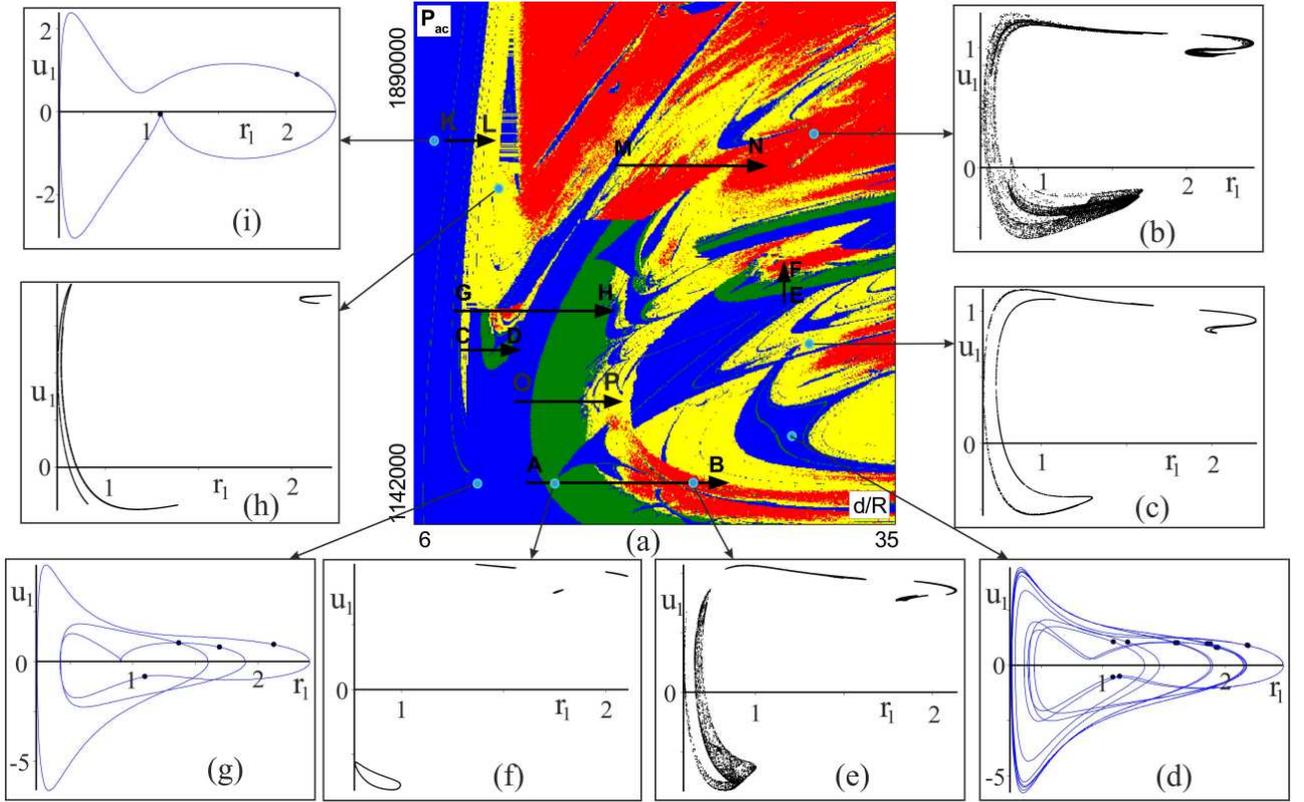} }
\caption{{\footnotesize (a) Charts of Lyapunov exponents for $\omega = 2.87 \cdot 10^7$ s$^{-1}$ and (b)--(i) projections of phase portraits of steady state regimes (attractors) for some representative points of this chart. Projections of attractors on Poincar\'e map are pained in black color, projections of phase trajectories for periodic attractors -- in blue color. The following attractors are shown here: (b) hyperchaotic attractor at $d/R_0 = 32, P_{ac} =1.68$ MPa with largest Lyapunov exponents of $\lambda_1 = 0.0803, \ \lambda_2 = 0.0357$; (c) synchronous chaotic attractor at $d/R_0 = 30, P_{ac} =1.4$ MPa with $\lambda_1 = 0.0684, \ \lambda_2 = -0.0268$; (d) syncronous 12-periodic limit cycle at $d/R_0 = 28, P_{ac} =1.3$ MPa with $\lambda_1 =-0.0616, \ \lambda_2 = -0.0733$; (e)  hyperchaotic attractor at $d/R_0 = 22, P_{ac} =1.2$ MPa with $\lambda_1 = 0.0241, \ \lambda_2 = 0.0034$; (f) quasiperiodic attractor at $d/R_0 = 14.5, P_{ac} =1.2$ MPa with $\lambda_1 = 0, \ \lambda_2 = -0.0149$; (g)  asynchronous 4-periodic limit cycle at $d/R_0 = 10, P_{ac} =1.2$ MPa with $\lambda_1 = -0.2331, \ \lambda_2 = -0.2343$; (h) synchronous chaotic attractor at $d/R_0 = 10, P_{ac} = 1.6$ MPa with $\lambda_1 = 0.0802, \ \lambda_2 = -0.0826$; (i) synchronous 2-periodic limit cycle at $d/R_0 = 6.75, P_{ac} = 1.7$ MPa with $\lambda_1 = -0.1437, \ \lambda_2 = -0.2057$.
}}
\label{fig:clock}
\end{figure}

Fig.~\ref{fig:clock}a shows a chart of two maximal Lyapunov exponents $\lambda_1 \geq \lambda_2$ on $(d/R_0, P_{ac})$ parameter plane for fixed $\omega = 2.87 \cdot 10^7$ s$^{-1}$. This value of $\omega$ belongs to the range of frequencies which is relevant for biomedical applications. Also, the system demonstrates quite rich dynamics at this frequency, and typical attractors and bifurcation scenarios are presented in physically relevant interval of pressures and distances between the bubbles.

Depending on values of $\lambda_1$ and $\lambda_2$ the corresponding pixel on the chart is painted with a certain color using the following scheme:
\begin{itemize}
\item $\lambda_1 < 0, \lambda_2 < 0$ -- periodic regime -- blue color;
\item $\lambda_1 = 0, \lambda_2 < 0$ -- quasiperiodic regime -- green color;
\item $\lambda_1 > 0, \lambda_2 \leq 0$ -- simple chaotic regime (strange attractor with one positive Lyapunov exponent) -- yellow color;
\item $\lambda_1 > 0, \lambda_2 > 0$ -- hyperchaotic regime (strange attractor with two positive Lyapunov exponents) -- red color;
\end{itemize}

Big blue region on the left in the chart of Lyapunov exponents corresponds to the stable periodic regimes (see Fig.-s~\ref{fig:clock}d, i), when two bubbles perform in-phase synchronous oscillations. At the top part of the chart this blue region adjoins yellow domain corresponding to a simple strange attractor with one positive Lyapunov exponent  (see Fig.~\ref{fig:clock}h), corresponding to synchronous chaotic oscillation of two bubbles.  In Sec.~\ref{sec:4.1.1} we show that the transition to chaos in this part of the chart occurs via Feigenbaum's cascade of period-doubling bifurcations.

At the middle and bottom parts of the diagram from Fig.~\ref{fig:clock}a the big blue region adjoins the green-colored region corresponding to quasiperiodic regime, when torus (invariant curve on the Poincar\'e map) is an attractor in the system, see Fig.~\ref{fig:clock}f. The invariant torus in Fig. \ref{fig:clock}f is born from the asynchronous periodic regime shown in Fig. \ref{fig:clock}g via the Neimark--Sacker bifurcation. It is worth noting that the asynchronous periodic limit cycle in Fig. \ref{fig:clock}g is (likely) born via a saddle-node bifurcation, which is not reflected in the two-dimensional chart. Thus if we move left (decrease $d/R_0$) from the point corresponding to Fig. \ref{fig:clock}g to the border of the plot, there will be a switch from asynchronous limit cycle to a synchronous one (the point where the saddle-node bifurcation occurs), and this switch cannot be observed in the two-dimensional chart  Fig.~\ref{fig:clock}a. Moving through the region corresponding to the quasiperiodic regime from left to right one can observe transition to chaotic (in the middle part of the chart) and even hyperchaotic (in the bottom part of the cart) attractors.
Simple chaotic attractors occurring after the destruction of an invariant torus as well as hyperchaotic attractors, correspond to asynchronous oscillations of the bubbles, see Fig.~\ref{fig:AS}c. In Sec.~\ref{sec:4.1.2} we show that in the middle part of the chart chaotic attractors appear in accordance with Afraimovich--Shilnikov scenario of the destruction of an invariant torus. In Sec.~\ref{sec:hyperchaos} we study the scenario for the onset of hyperchaotic oscillations in the bottom part of the chart.

From Fig.~\ref{fig:clock}a one can see that yellow- and red-colored regions corresponding to chaotic (see e.g. Fig.~\ref{fig:clock}c,h) and hyperchaotic (see e.g. Fig.~\ref{fig:clock}b,e) attractors alternate with the so-called stability windows inside which stable periodic orbits are observed (see Fig.~\ref{fig:clock}a). Some of these stability windows have shrimp-like form \cite{Gallas1993}. Such stability windows indicate the existence of the specific homoclinic bifurcations (cubic homoclinic tangencies or symmetrical pairs of homoclinic tangencies) in the system \cite{GonSimoVieiro2013}.
All these stability windows indicate that chaotic dynamics in system \eqref{eq:eq1} are not hyperbolic and even not pseudohyperbolic \cite{TurShil98, GonGonKazKoz18, GonKazTur18}. In other words, in the accordance with PQ-hypothesis \cite{GonKazTur18} strange attractors in the system under 
investigation belong to a class of quasiattractors introduced by Afraimovich and Shilnikov in \cite{AfrSh83a}. Stable periodic orbits of high periods and with narrow absorbing domains exist inside such attractors or appear with arbitrarily small perturbations. However, from a physical point of view in most cases quasiattractors do not differ from genuine (pseudohyperbolic) attractors due to narrow absorbing domains of periodic orbits belonging to them. It is important to note that currently there are no known systems which demonstrate hyperbolic or even pseudohyperbolic hyperchaotic behavior.


\section{Transition to chaos and hyperchaos} \label{sec:4}

As it can be clearly seen from the chart of Lyapunov exponents (Fig.~\ref{fig:clock}a) the dynamics in system \eqref{eq:eq1} can be either regular (periodic or quasiperiodic) or chaotic and even hyperchaotic.

As a rule, chaotic attractors appear from simple (regular) attractors as a result of the implementation of some bifurcation scenario. The most known examples of such scenarios are: 1) the Feigenbaum's cascade of period-doubling bifurcations \cite{Feigenbaum} according to which chaotic attractor appears from a stable periodic orbit via infinite sequence of period doubling bifurcations; 2) destruction of an invariant torus by Afraimovich-Shilnikov scenario \cite{AfrShil1983}; 3)  the Shilnikov scenario \cite{Shil86} due to which spiral attractor containing a saddle-focus equilibrium with a two-dimensional unstable invariant manifold appears from a stable equilibrium as a result of certain local and global bifurcations. It is worth noting that all the above mentioned scenarios can be observed in flow system with dimension $N \geq 3$ or in diffeomorphisms -- discrete system (except for the Shilnikov scenario) with dimension $N \geq 2$ and generally lead to the appearance of chaotic attractors with only one positive Lyapunov exponent.

An important class of chaotic attractors of multi-dimensional flows ($N \geq 4$) and maps ($N \geq 3$), namely so-called \textit{homoclinic attractors} containing saddle periodic orbits with its homoclinic structure, was introduced in \cite{GonGonShil2012}, \cite{GonGonKazTur2014}, where the classification of such attractors and also phenomenological scenarios of their appearance were proposed (see, also \cite{GonGonPhD16} for more examples of such attractors in three-dimensional H\'enon maps). The classification of homoclinic attractors is based on the type of a saddle orbit belonging to an attractor. Two main attractors of this type are the \textit{discrete Lorenz} and \textit{figure-eight} attractors. They contain a saddle fixed point with a one-dimensional unstable manifold forming a homoclinic structure resembling a butterfly and figure-eight, respectively. As it was recently shown in \cite{GonGonKazKoz18} these two attractors belong to a class of \textit{pseudohyperbolic} \cite{TurShil98}, \cite{ShilTur2008} (``genuinely'' chaotic) attractors. Shortly speaking, each orbit on a pseudohyperbolic attractor has a positive Lyapunov exponent and, what is important from a physical point of view, this property persists after small perturbations (changing in parameters). However, both the discrete Lorenz and figure-eight attractors cannot be hyperchaotic.

Another important example of a homoclinic strange attractor proposed in \cite{GonGonShil2012}, \cite{GonGonKazTur2014} is a \textit{discrete Shilnikov attractor}. In contrast to the all above mentioned examples of strange attractors, the discrete Shilnikov attractor contains a saddle-focus fixed point with two-dimensional unstable invariant manifold. In all small neighborhoods of such fixed point, two-dimensional areas are expanded. Lyapunov exponents on this attractor as a whole 'can feel' this expansion, and thus, two Lyapunov exponents can be positive. Unlike the discrete Lorenz and figure-eight attractors, the discrete Shilnikov attractor cannot be pseudohyperbolic. This attractor belongs to another class of strange attractors called by Afraimovich and Shilnikov in \cite{AfrSh83a} as \textit{quasiattractors}. Homoclinic tangencies inevitably arise in quasiattractors and lead to the birth of stable periodic orbits inside such attractors. Thus, these attractors either contain a stable periodic orbit with large periods and narrow absorbing domains or such orbits appear after arbitrary small perturbations (parameter changing).

Discrete Shilinikov attractors were found in different dynamical systems such as the generalized three-dimensional H\'enon maps, nonholonomic models of Chaplygin top \cite{BorKazSat2016} and Celtic stone \cite{GonGonKazSam18}, the model of coupled identical oscillators \cite{GrinKazSat2017}.
 However, not in all cases they were identified as hyperchaotic attractors. Apparently, in some cases the expansion of two-dimensional areas near a saddle-focus orbit with a two-dimensional unstable manifold is compensated by the volume contraction near some other saddle periodic orbits that also belong to the attractor, but which have only a one-dimensional unstable manifold. Since Lyapunov exponents are average characteristic of an attractor, only one Lyapunov exponent can become positive in this case.

In Sec.~\ref{sec:hyperchaos} we show that for system \eqref{eq:eq1} discrete Shilnikov attractors containing a saddle-focus periodic orbit are hyperchaotic. We also propose a new phenomenological scenario which leads to the appearance of hyperchaotic attractor and demonstrate that exactly due to this scenario hyperchaos appears in the system under consideration. However first of all we describe scenarios of transition to simple chaotic (with only one positive Lyapunov exponent) attractors in the model.

\subsection{Transition to chaos in the model} \label{sec:chaos}

\subsubsection{Feigenbaum's cascades of period doubling bifurcations} \label{sec:4.1.1}

Feigenbaum's infinite sequence of period-doubling bifurcations \cite{Feigenbaum} is one of the typical scenarios of the chaos onset in one- and two-dimensional maps and three-dimensional flows. However, such scenario is also found in multi-dimensional maps ($N \geq 3$) and flows ($N \geq 4$) (see, e.g. \cite{BorKazKuz14}). On the other hand, since in multi-dimensional systems period-doubling bifurcations compete with Neimark-Sacker bifurcations, the transition to chaos via Feigenbaum's cascade for multi-dimensional systems is more rare than in the case of two-dimensional maps and three-dimensional flows.

\begin{figure}[!ht]
\center{\includegraphics[width=1.0\linewidth]{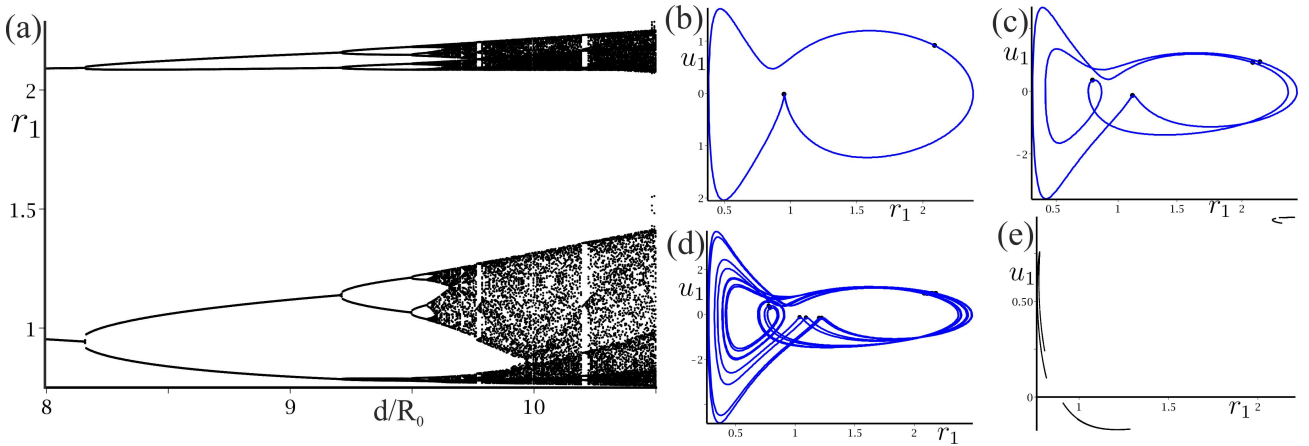} }
\caption{{\footnotesize Transition from synchronous periodic to synchronous chaotic oscillations on the path: $P_{ac} = 1.7$ MPa, $8.0 < d/R_0 < 10.5$ via the Fegeinbaum`s cascade. (a)  bifurcation tree; (b)  projection of the phase portrait of the 2-periodic limit cycle at $d/R_0 = 8$; (c)  4-periodic limit cycle at $d/R_0 = 9$; (d) 16-periodic limit cycle at $d/R_0 = 9.55$; (e)  projection of the Poincar\'e section of the chaotic attractor at $d = 9.85 \cdot R_0$ with two largest Lyapunov exponents: $\lambda_1 = 0.0501, \lambda_2 = -0.1339$.
}}
\label{fig:scenario_1}
\end{figure}

Here we show that strange attractors in system \eqref{eq:eq1} can appear via Feigenbaum's cascade of period-doubling bifurcations. Let us fix $P_{ac} = 1.7$ and move along the path KL
: $(P_{ac}, d/R_0) = (1.7 \cdot 10^6, 8) \rightarrow (P_{ac}, d/R_0) = (1.7 \cdot 10^6, 10.5) $. We start from $d/R_0 = 8$ because periodic attractors at $d/R_0=6.75$ (see Fig. \ref{fig:clock}i) and $d/R_0 = 8$ (see Fig. \ref{fig:scenario_1}b) look identical and we do not lose any information by starting from this point. Fig.~\ref{fig:scenario_1}a shows the corresponding bifurcation tree in this case. Phase portraits of attractors at some values of parameter $d/R_0$ are presented in Figs.~\ref{fig:scenario_1}b--d, where in blue color we show the projection of the phase curves onto the $(r_1, u_1)$ plane, and in black color -- projections of the corresponding Poincar\'e map (at $t=2 \pi k, k \in \mathbb{N}$) on the same plane. 
At the starting point of the path a stable periodic orbit (point of period 2 on the Poincar\'e map) is an attractor of the system, see Fig.~\ref{fig:scenario_1}b. When the parameter $d/R_0$ increases this periodic orbit undergoes cascade of period doubling bifurcations, see Figs.~\ref{fig:scenario_1}c,d and finally, at $d/R_0 \approx 9.6$, Feigenbaum-like strange attractor emerges, see Fig.~\ref{fig:scenario_1}e where the corresponding Poincar\'e map is presented. The set of Lyapunov exponents for this attractor at $(d/R, P_{ac}) = (9.85, 1.7)$ is $\lambda_1 = 0.0501, \lambda_2 = -0.1339, \lambda_3 =-0.3733, \lambda_4 = -0.4881$.

In the general case Feigenbaum's cascade gives rise to the onset of a strange attractor with only one positive Lyapunov exponent. However, it is important to note, that at some specific cases such scenario can lead to the onset of hyperchaotic behavior. For example, if we take two identical oscillators demonstrating transition to chaos via cascade of period-doubling bifurcations and make them interact through a very weak coupling, both elements will demonstrate chaotic behavior, which will lead to two positive Lyapunov exponents in the coupled system. Such transition to hyperchaos was observed e.g. in \cite{Stankevich}. Since we take two identical elements the transition to hyperchaos in accordance with this scenario is also possible in our system if we suppose that the bubbles are quite distant (i.e. $d/R_{0}\ll 1$). However, this case is less interesting from a physical point of view and, therefore, we do not consider it here.

\subsubsection{Afraimovich-Shilnikov scenario of torus destruction} \label{sec:4.1.2}

As one can see from the chart of Lyapunov exponents (Fig. \ref{fig:clock}) some regions of parameters corresponding to stable periodic regimes adjoint the region with stable qusiperiodic regimes -- invariant tori. The boundary between these regions is formed by the curve of supercritical Neimark-Sacker bifurcation. Passing through this curve stable limit cycle loses its stability, becomes of a saddle-focus type with a two-dimensional unstable invariant manifold and a stable invariant torus appears.

From another side of the regions of the existence of quasiperiodic regimes the dynamics of system \eqref{eq:eq1} can be chaotic. This means that in system \eqref{eq:eq1} chaotic attractors can appear after destruction of a torus. There are a few typical scenarios of transition to chaos through the destruction of an invariant torus. One of such scenarios was proposed by Afraimovich and Shilnikov in \cite{AfrShil1983}.

\begin{figure}[!ht]
\center{\includegraphics[width=1.0\linewidth]{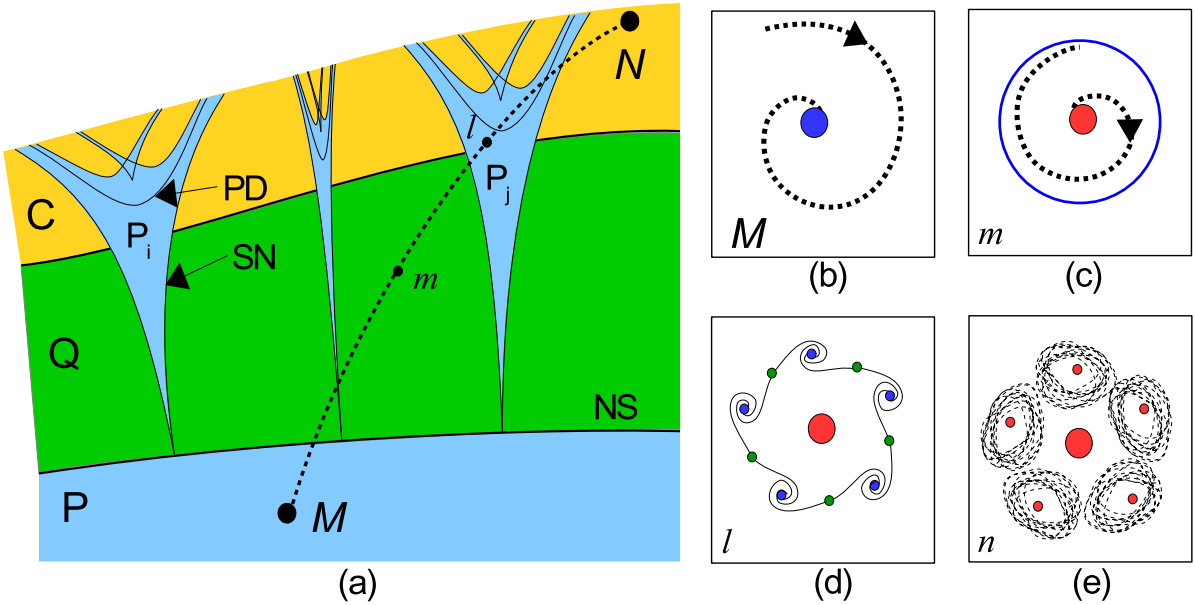} }
\caption{{\footnotesize (a) Sketch of the bifurcation diagram illustrating bifurcation of an invariant torus. $P$, $Q$, $P_i$ and $C$ -- regions of the existence of (b) stable periodic orbit, (c) stable invariant torus, (d) resonant periodic orbits, and (e) torus-chaos attractors, respectively. $MN$ -- some path along which torus-chaos attractor appears in accordance with Afraimovich-Shiknikov scenario.
}}
\label{fig:torus_2d}
\end{figure}

Here we show that the Afraimovich--Shilnikov scenario is the second typical scenario (the first is Feigenabum's cascade of period-doubling bifurcations) for the onset of chaos in system \eqref{eq:eq1}.
But first of all, let us recall some important details concerning typical organization of a bifurcation diagram inside region $Q$ of an invariant torus existence for two-dimensional maps, see Fig.~\ref{fig:torus_2d}a. Resonance regions $P_i$ -- the so-called tongues of synchronization alternate with quasiperiodic regions $Q$ above the curve of the Neimark-Sacker bifurcation $NS$ and with chaotic regions in the upper part of the diagram. Notice that in the tongues of synchronization resonant stable and saddle periodic orbits appear on torus (through the saddle-node bifurcations $SN$), while this torus still exists, but now it is formed by the closure of the unstable invariant manifold of the resonant saddle orbit\footnote{For multi-dimensional maps ($N \geq 3$) the structure of bifurcation diagrams is similar but slightly different from two-dimensional case. Inside tongues of synchronization, together with period-doubling bifurcations, Neimark-Sacker bifurcations are also possible, see Fig.~\ref{fig:torus_3d}a.}, see Fig.~\ref{fig:torus_2d}d.
\begin{figure}[!hb]
\center{\includegraphics[width=1.0\linewidth]{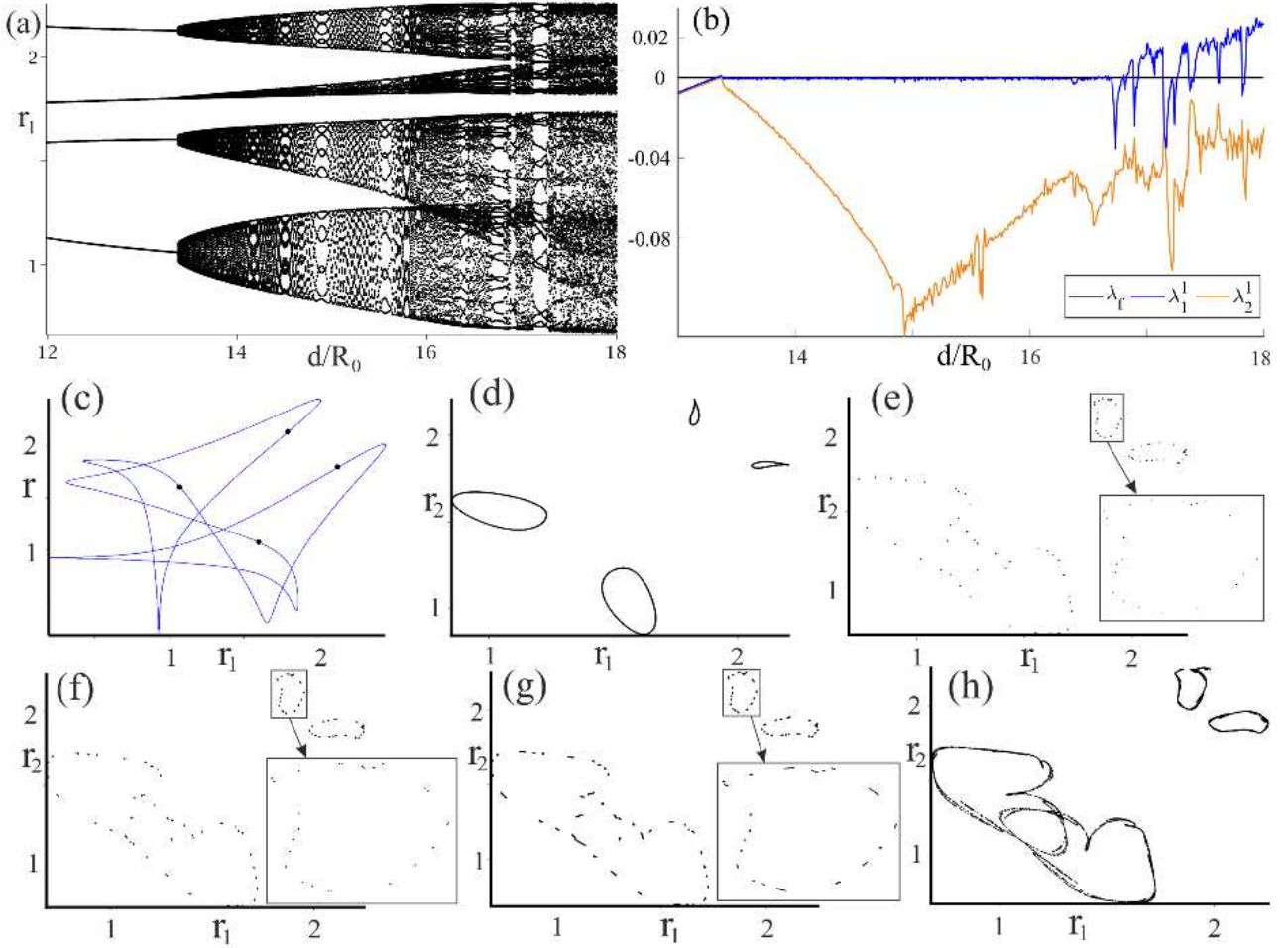} }
\caption{{\footnotesize Afraimovich--Shilnikov scenario for the onset of torus-chaos attractor along the path OP: $P_{ac} =1.35$ MPa, $d/ R_0 \in [13, 18]$. (a) bifurcation tree and (b) map of the two largest Lyapunov exponents associated with this path; (c)--(h) Projections of attractors onto $(r_1, r_2)$-plane: (c)  stable limit cycle (point of period 4 on the Poincar\'e map) at $d/R_0 = 13$; (d)  four-component torus after the Neimark-Sacker bifurcation at $ d/R_0 = 14.5$; (e)   resonance on the torus,  $d/R_0 = 16.825$; (f)  resonance after the first period-doubling bifurcation, $d/R_0 = 16.831$; (g)  $d/R_0 = 16.84$, 
Feigenbaum's cascade continunues leading to chaos ; (h) --  torus-chaos at $d/R_0 = 17.02$ with the following Lyapunov spectre $\lambda_1 = 0.0139, \lambda_2 = -0.0538, \lambda_3 = -0.5571, \lambda_4 = -0.6131$.
}}
\label{fig:AS}
\end{figure}
Moving along arbitrary path in the parameter plane one can observe sequences of alternated regular, quasiperiodic and chaotic regimes. Thus, it is important to note that bifurcations of an invariant torus and, in particular, transition to chaos depend on the path in the bifurcation diagram. Moreover, the parts of this path from resonance regions to chaotic ones are the most important. For example, moving along path $MN$ one can observe the following sequence of regimes: stable periodic orbit (Fig.~\ref{fig:torus_2d}b), stable torus (Fig.~\ref{fig:torus_2d}c), resonance torus existing in resonance region $P_j$ (Fig.~\ref{fig:torus_2d}d), and, finally, torus-chaos attractor (Fig.~\ref{fig:torus_2d}e).


According to Afraimovich and Shilnikov \cite{AfrShil1983}, in two-parametric families of two-dimensional maps, destruction of an invariant torus inside resonance regions can happen due to the following scenarios: 1) period-doubling bifurcation with a stable resonant orbit (e.g. if to move upwards inside a resonance region); 2) homoclinic bifurcation, when unstable invariant manifold of the resonant saddle periodic orbit touches (and than intersects) its stable manifold (e.g. if to move inside a resonance region through the curve $SN$ to the chaotic region, see path $MN$ in Fig.~\ref{fig:torus_2d}a); 3) more complex and difficult to observe scenario associated with the increasing of oscillations of the unstable manifold of a resonant saddle orbit, see details in \cite{AfrShil1983}.

What is important, in all these cases, when leaving a resonance region, one can observe the chaotic regime associated with the previously existed invariant torus. Such chaotic attractors were called as \textit{torus-chaos attractor} in \cite{AfrShil1983}.

Fig.~\ref{fig:AS} gives an illustrative example of the onset of torus-chaos attractor in system \eqref{eq:eq1} in accordance with Afraimovich-Shilnikov scenario. Here we fix $P_{ac}= 1.35$ MPa and move along path OP from the chart of Lyapunov exponents (Fig.~\ref{fig:clock}): $13 < d/R_0 < 18$. Fig.-s~\ref{fig:AS}a and \ref{fig:AS}b show the corresponding bifurcation tree and the graph of the two largest Lyapunov exponents in this case. Portraits of some attractors along this path are presented in Fig.-s~\ref{fig:AS}c--h, where in blue color we show the projection of phase portraits onto $(r_1, r_2)$ plane, and in black color -- projections of the corresponding Poincar\'e map on the same plane $(r_1, r_2)$.

The stable limit cycle (stable fixed point of period four on the Poincar\'e map) existing at $d/R = 13$
 (see Fig.~\ref{fig:AS}c) undergoes the Neimark-Sacker bifurcation at $d/R \approx 13.28$ after which 
 a stable invariant torus (four-component invariant curve in the Poincar\'e map) appears, see Fig.~\ref{fig:AS}d. Then, at $d/R \approx 16.71$ we get into a resonance region where stable periodic orbit appears, see Fig.~\ref{fig:AS}e. Moving inside this resonance region the stable periodic orbit undergoes Feigebaum's cascade of period-doubling bifurcations, see Fig.-s~\ref{fig:AS}e--g for details.
 Finally, at $d/R \approx 16.84$, we get out of the resonance region and torus-chaos attractor appears, see Fig.~\ref{fig:AS}g.

In the next subsection we will demonstrate that some paths out of resonance regions lead to the onset of hyperchaotic attractors. We also will give a bifurcation scenario of such transition.

\subsection{Transition to hyperchaos in system \eqref{eq:eq1}} \label{sec:hyperchaos}

Starting from three-dimensional maps (four-dimensional flows), in addition to Afraimovich-Shilnikov scenarios, some other ways of the destruction of an invariant torus become possible. Here we would like to mention the following two scenarios: 1) cascades (finite) of period-doubling bifurcation of an invariant torus \cite{ArneodoCoulletSpiegel83}, \cite{Kaneko}, \cite{Anishchenko05}; 2) secondary Neimark-Sacker bifurcation with a stable resonant periodic orbit inside a tongue of synchronization.

\begin{figure}[!ht]
\center{\includegraphics[width=1.0\linewidth]{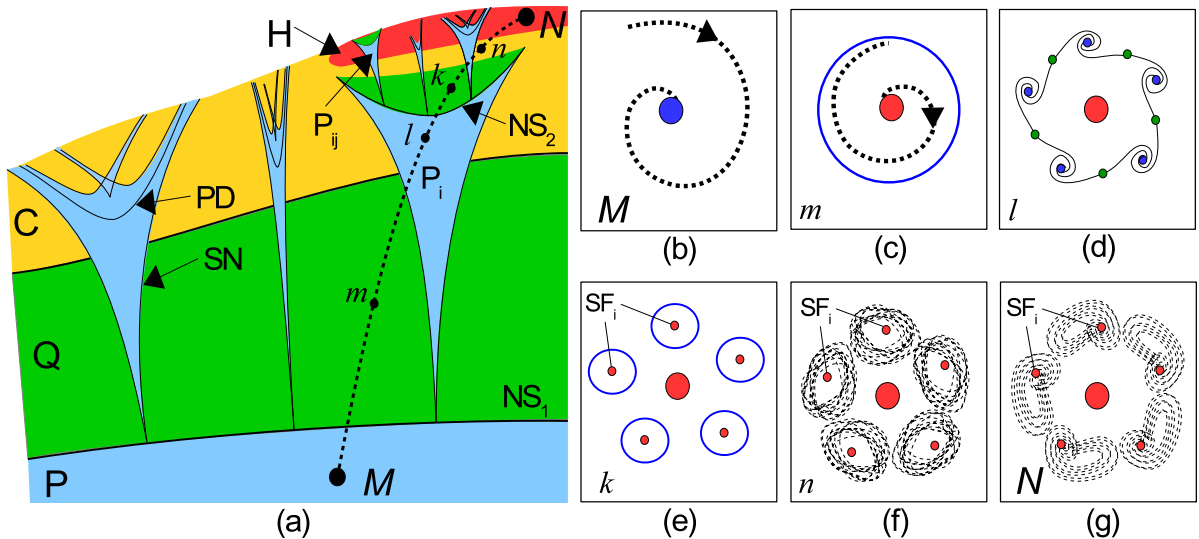} }
\caption{{\footnotesize Sketch of the bifurcation diagram illustrating the scenario of the onset of a hyperchaotic attractor in multi-dimensional maps ($N \geq 3$). $P$, $Q$, $P_i$, $C$ and $H$ -- regions of the existence of (b) stable periodic orbit, (c) stable invariant torus, (d) resonant periodic orbits, (e) stable torus after secondary Neimark-Sacker bifurcation $NS_2$, (f) torus-chaos attractors, and (g) hyperchaotic Shilnikov attractor, respectively. $MN$ -- some path along which hyperchaotic attractor appears.
}}
\label{fig:torus_3d}
\end{figure}

Thus, a typical bifurcation diagram near the Neimark-Sacker bifurcation for three-dimensional maps differ from the corresponding diagram in two-dimensional case, see Fig.~\ref{fig:torus_3d}a. The main difference is that a resonant periodic orbit in three-dimensional case can undergo the secondary Neimark-Sacker bifurcation $NS_2$ instead of a typical for two-dimensional maps period-doubling bifurcation, see right-top part in Fig.~\ref{fig:torus_3d}a.

Secondary Neimark-Sacker bifurcation is the first (but not main) step in our scenario for the onset of hyperchaotic attractors. After this bifurcation multi-round stable invariant torus (multi-component invariant curve in the Poincar\'e map) is born in the system, while, what is very important, resonant periodic orbit becomes a saddle-focal with a two-dimensional unstable manifold, see Fig.~\ref{fig:torus_3d}e, where the saddle-focus periodic orbit is denoted as $SF_i$.

The next step in the framework of this scenario is associated with the destruction of a multi-round stable invariant torus. It does not matter how it happens through the Afraimovich-Shilnikov scenario, cascade of torus period-doubling bifurcations or even via the tertiary Neimark-Sacker bifurcation. But we suppose, and it is quite natural, that after the corresponding bifurcations torus-chaos attractor with one positive Lyapunov exponent appears, see Fig.~\ref{fig:torus_3d}f. It is worth noting that in this case, immediately after transition to chaos, saddle-focus $SF_i$ does not belong to the torus-chaos attractor. It means that orbits of this chaotic attractor do not attend some neighborhood of $SF_i$, see Fig.~\ref{fig:torus_3d}f.


The final, and the key step in the scenario is the inclusion of the saddle-focus periodic orbit $SF_i$ which appeared after the secondary Neimark-Sacker into the chaotic attractor. After this inclusion, saddle-focus orbit $SF_i$ together with its two-dimensional unstable manifold and its homoclinic structure starts to belong to the attractor, i.e. discrete homoclinic Shilnikov attractor based on this saddle-focus orbit emerges, see Fig.~\ref{fig:torus_3d}g. Orbits on this attractor can pass arbitrary close to $SF_i$, where two-dimensional areas are expanded. As a result, two Lyapunov exponents become positive i.e. a hyperchaotic attractor is born.

The inclusion of a saddle-focus periodic orbit to the chaotic attractor can occur in different ways. It depends on the transition from the stable multi-round torus to chaotic regime. In all known models demonstrating the onset of the discrete Shilnikov attractor (in three-dimensional H\'enon maps \cite{GonGonPhD16}, nonholonomic models of Chaplygin top \cite{BorKazSat2016} and Celtic stone \cite{GonGonKazSam18}) this inclusion happens in a soft manner by a smooth transformation of a torus-chaos attractors. However, we also suppose that a saddle-focus orbit can be included to the chaotic attractor sharply
 due to the crisis of multi-round torus-chaos attractor.

\begin{figure}[!ht]
\center{\includegraphics[width=1.0\linewidth]{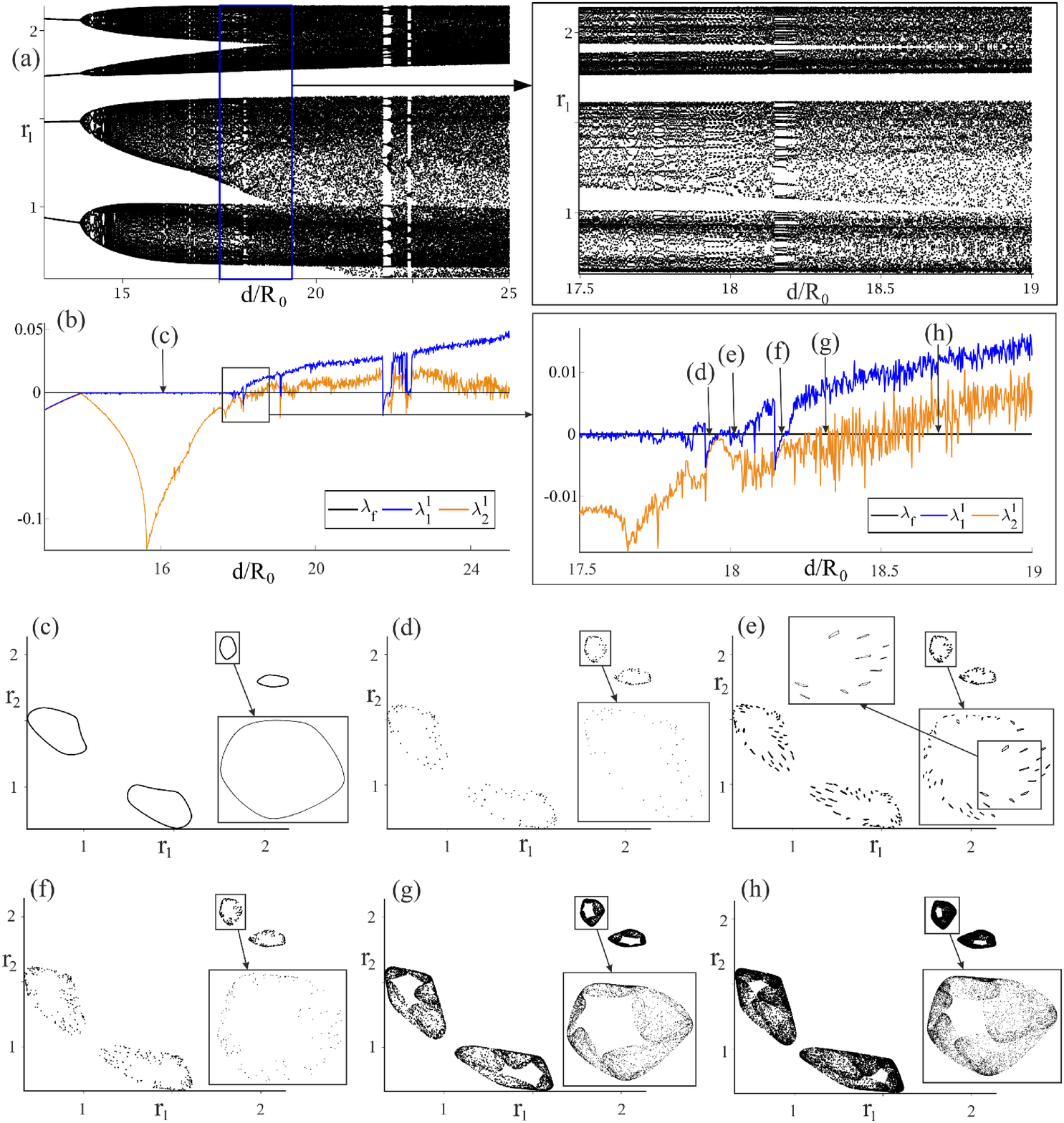} }
\caption{{\footnotesize The implementation of the proposed scenario of the onset of a hyperchaotic attractor along the path $AB$: $P_{ac} = 1.2$, $13 < d/R_0 < 25$. (a)  and (b) bifurcation tree and the graph of two largest Lyapunov exponents associated with this path with the enlarged area for $17.5 < d/R_0 < 19$; (c)--(h) projections of the Poincar\'e maps for different attractors on the $(r_1, u_1)$ plane: (c) four-component torus at $d/R_0 = 16$; (d)  resonance orbit at $d/R_0 = 17.92$; (e)  multi-component torus after the secondary Neimark-Sacker bifurcation, $d/R_0 = 17.98$; (f)  high-period resonance orbit emerges on the secondary torus and starts to go through the period-doubling bifurcations, $d/R_0 = 18.01$; (g) -- chaotic attractor after the period-doubling cascade at $d/R_0 = 18.35$ with two largest Lyapunov exponents $\lambda_1 = 0.0087 , \lambda_2 = -0.0040$, one can see the gaps in the chaotic attractor existing around the saddle-focus orbit 
; (h) hyperchaotic Shilnikov attractor containing saddle-focus periodic orbit with a two-dimensional unstable manifold at $d/R_0 = 18.71$.
}}
\label{fig5a}
\end{figure}


In order to support the proposed scenario we show that the transition to hyperchaos along paths $AB$ and $GH$ (see Fig. \ref{fig:clock}a) happens in full compliance with this scenario. Firstly, let us consider the path $AB$, corresponding to the following parameters interval: $P_{ac} = 1.2$ MPa and $13 < d/R_0 < 25$. The bifurcation tree, corresponding to this route is shown in Fig. \ref{fig5a}a and the graph of two largest Lyapunov exponents is presented in Fig. \ref{fig5a}b. We also show the enlarged part for both these graphs at the right panels in Fig.-s \ref{fig5a}a,b in order to explore some important details concerning secondary Neimark-Sacker bifurcation and the transition to hyperchaos. Projections of the Poincar\`e sections for some representative attractors are shown in Fig.-s \ref{fig5a}c--h.
From Fig. \ref{fig5a} one can observe the following bifurcations sequence happening to the asynchronous 4-periodic limit cycle existing at $d/R_0 = 13$:  
 it  undergoes the Neimark-Sacker bifurcation at $d/R_0 \approx 13.89$ and the 4-component torus arises, see Fig. \ref{fig5a}c. Then a high-periodic resonance occurs on the torus, see Fig. \ref{fig5a}d. With further increasing of $d/R_0$ the multi-component torus emerges after the secondary Neimark-Sacker bifurcation while the former resonance orbit becomes saddle-focus with two-dimensional unstable manifold, see \ref{fig5a}e. Soon the multi-component torus gives rise to the torus-chaos attractor (see Fig. \ref{fig5a}g), which appears via the cascade of period-doubling bifurcations happening with some stable resonant orbit emerging on this torus (see the long-periodic orbit emerging after few period-doubling bifurcations in Fig. \ref{fig5a}f). It is important to note that saddle-focus orbit occurring after the secondary Neimark-Sacker bifurcation does not belong to this torus-chaos attractor. Gazing at Fig. \ref{fig5a}g, one can see some visible gaps existing around the saddle-focus orbit.
 Finally, at $d/R_0 \approx 18.66$, the saddle-focus orbit starts to belong to the chaotic attractor. As a result the hyperchaotic attractor, containing this saddle-focus orbit appears, see Fig.\ref{fig5a}h. The Lyapunov exponents at $d = 18.71$ are $\lambda_1 = 0.0135 , \lambda_2 = 0.0019, \lambda_3 = -0.5560 , \lambda_4 = -0.5607$.

\begin{figure}[!ht]
\center{\includegraphics[width=1.0\linewidth]{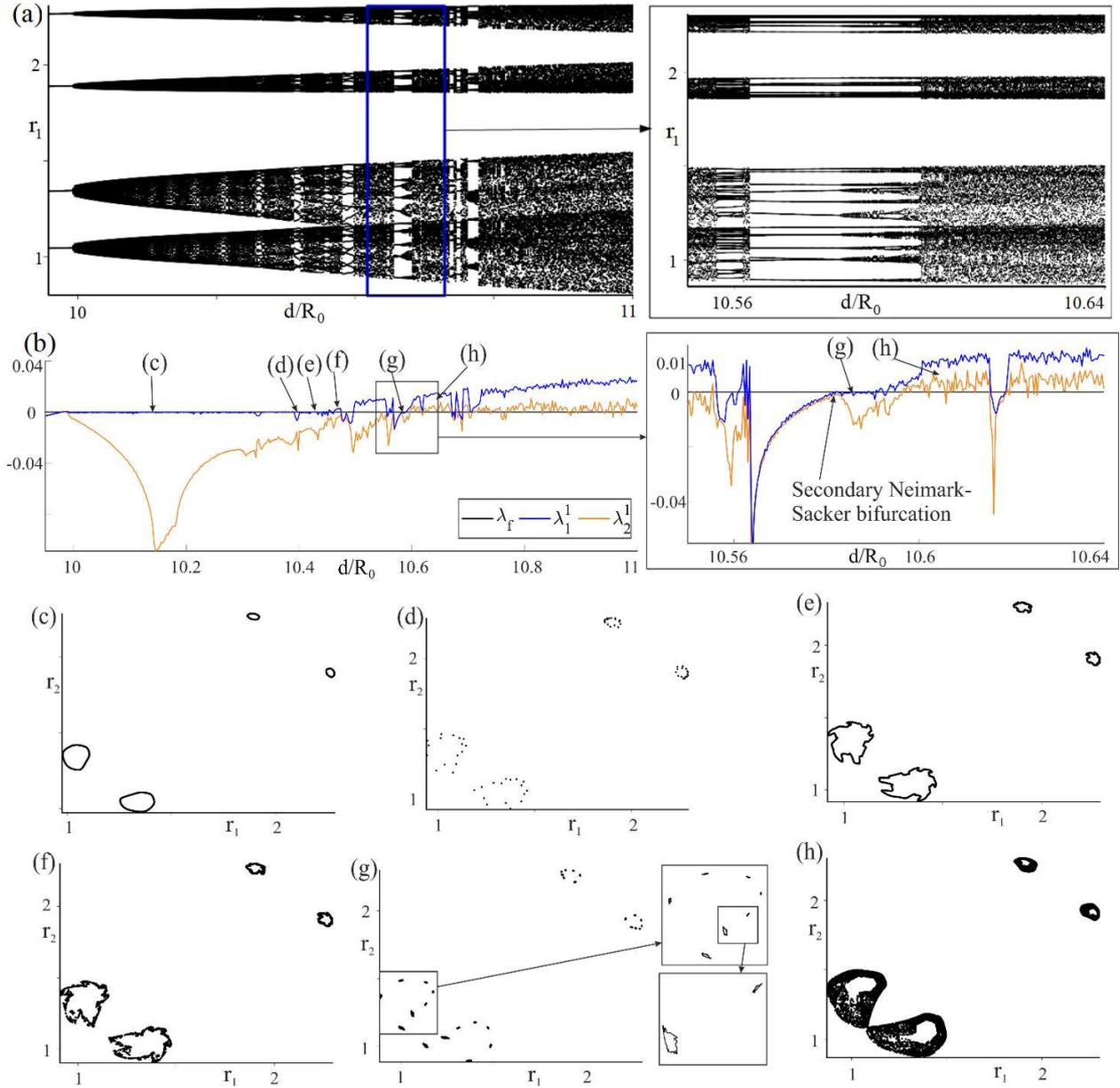} }
\caption{{\footnotesize The implementation of the same scenario of the onset of hyperchaotic attractor along the path $GH$: $P_{ac} = 1.45$ MPa, $9.85 < d/R_0 < 11$. (a) and (b) bifurcation tree and graph of two largest Lyapunov exponents associated with this path; (c)--(h)  projections of  the Poincar\'e map of several attractors on the $(r_1, r_2)$ plane. (c)  four-component torus at $d/R_0 = 10.15$; (d) resonance orbit at $d/R_0 = 10.3961$; (e)  torus starts losing its smoothness, $d/R_0 = 10.4111$; (f) torus-chaos attractor at $d/R_0 = 10.4687$; (g) multi-component torus after the secondary Neimark-Sacker bifurcation, $d/R_0 = 10.59$; (h) hyperchaotic attractor containing saddle-focus periodic orbit with a two-dimensional unstable manifold at $d/R_0 = 10.6115$.
}}
\label{fig5}
\end{figure}


Further we demonstrate several more routes illustrating the same scenario of the onset of a hyperchaotic attractor. Let us fix $P_{ac} = 1.45$ and move along the path $GH$ from the chart of Lyapunov exponents in Fig.~\ref{fig:clock}a. Fig.-s~\ref{fig5}a,b show the corresponding bifurcation tree and a the graph of the largest Lyapunov exponents along this path. From these figures one can see that at $d/R_0 \approx 10.59$ two Lyapunov exponents become positive, i.e. hyperchaotic attractor appears.  Portraits of some represntative attractors for several values of parameter $d/R_0$ are presented in Fig.-s~\ref{fig5}c--h, where projections of the Poincar\'e map onto the $(r_1, r_2)$ plane are shown. The beginning of the route GH corresponds to an asynchronous
 4-periodic limit cycle. At $d/R_0 \approx 9.98$ it undergoes the Neimark-Sacker bifurcation, after which a stable invariant torus (four component invariant curve on the Poincar\'e map) appears, see Fig.~\ref{fig5}c. With increasing $d/R_0$ we pass through few resonance regions (see Fig.~\ref{fig5}d). Soon after this the invariant torus starts to loss its smoothness, see Fig.~\ref{fig5}e. Then, leaving one of the resonance regions, torus-chaos attractor with one positive Lyapunov exponent and the following spectre appears: $\lambda_1 = 0.0090, \lambda_2 = -0.0404, \lambda_3 = -0.5302, \lambda_4 = -0.5507$, see Fig.~\ref{fig5}f. With further increasing  $d/R_0$ we again pass through a resonance region, but now the resonant orbit undergoes the secondary Neimark-Sacker bifurcation at $d/R_0 \approx 10.58$ after which a multi-round invariant torus (multi-component invariant curve on the Poincar\'e map) appears, while the resonance orbit becomes saddle-focus $SF$ with a two-dimensional unstable manifold, see Fig.~\ref{fig5}g. Not long after this multi-round invariant torus gives rise to a 
chaotic attractor 
 and finally the discrete hyperchaotic Shilnikov attractor containing the saddle-focus orbit $SF$ appears
 (see Fig.~\ref{fig5}g). The set of Lyapunov exponents for this hyperchaotic attractor at $(d/R, P_{ac}) = (10.6115, 1.45)$ is $\lambda_1 = 0.012, \lambda_2 = 0.0022, \lambda_3 = -0.559, \lambda_4 = -0.569$.

\begin{figure}[!ht]
\center{\includegraphics[width=1.0\linewidth]{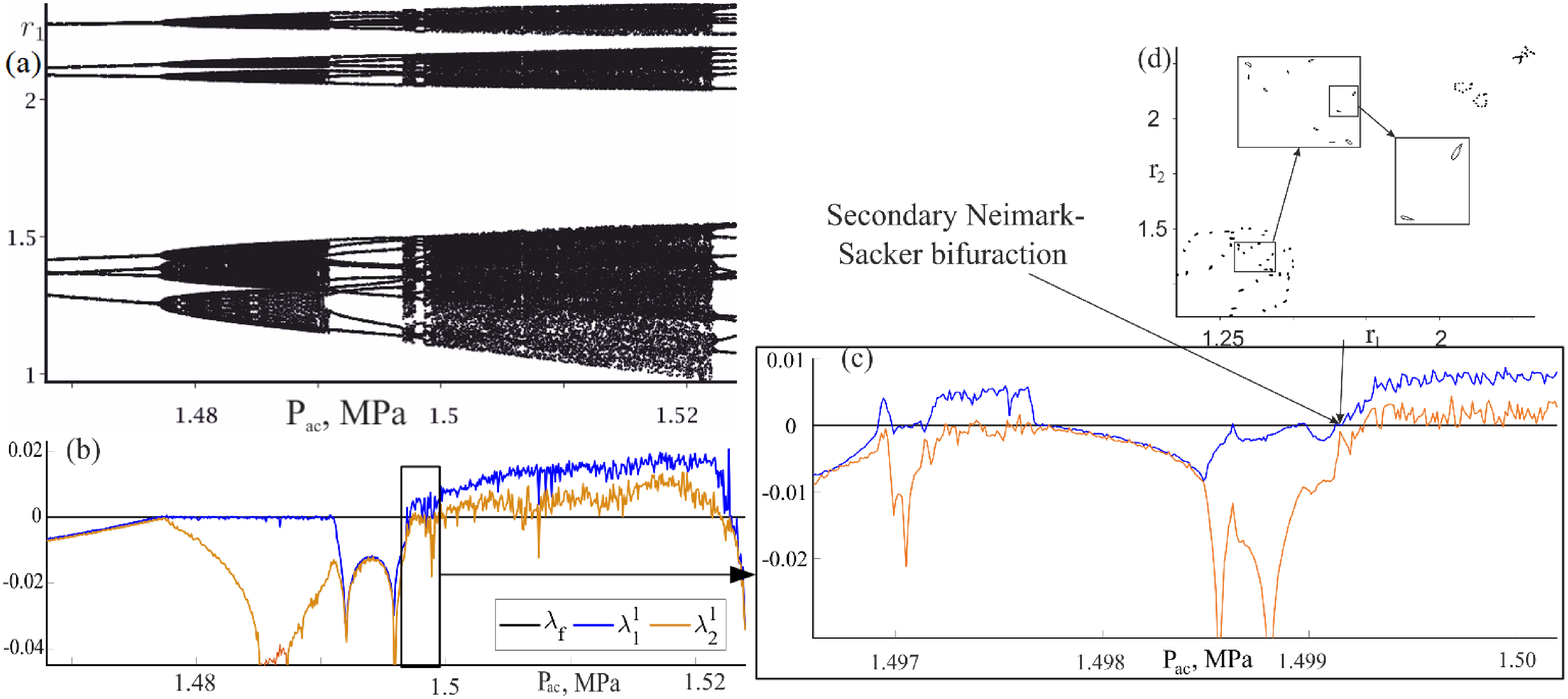} }
\caption{{\footnotesize Transition to hyperchaos along the route EF:  $d = 27.9 \cdot R_0$, $P_{ac}$: $1.468$ MPa$ < P_{ac} < 1.524$ MPa. (a), (b) bifurcation tree and graph of two largest Lyapunov exponents associated with this path; (c) enalarged part of the graph of the Lyapunov exponents, corresponding to $1.4966$ MPa $< P_{ac} < 1.5002$ MPa; (d) projection of the Poincar\'e map of the multicomponent torus after the secondary Neimark--Sacker bifurcation on the $(r_1, r_2)$ plane.
}}
\label{figEF}
\end{figure}


The third route we breifly discuss here is the route EF: $d/R_0 = 27.9, 1.468$ MPa $< P_{ac} < 1.524$ MPa. Unlike for the previous routes, for this one we fix $d/R_0$ and vary $P_{ac}$.
In Fig. \ref{figEF}a,b we present the bifurcation diagram and the graph of two largest Lyapunov exponents associated with the route. At the beginning of this path there is an asynchronous 6-periodic limit cycle. Increasing $P_{ac}$ we can observe the first Neimark-Sacker bifurcation occuring at $P_{ac} \approx 1.477$ MPa (see Fig. \ref{figEF}b). The bifurcation scenario proposed earlier is very hard to notice, unless we study a significantly enlarged picture, see \ref{figEF}c. Further increasing $P_{ac}$ leads to emergence of a resonance. Soon, at $P_{ac} = 1.49914$ MPa it undergoes the secondary Neimark-Sacker bifurcation at which the former resonance orbit becomes saddle-focus with two-dimensional unstable manifold and a multi-component torus emerges, see \ref{figEF}d. This secondary torus is very hard to notice because it exists in a very narrow range of pressures. Further increasing of $P_{ac}$ leads to inclusion of the saddle-focus orbit with two-dimensional unstable manifold into the attractor and it becomes hyperchaotic.


Finally, we would like to raise an open question. In which cases a discrete Shilnikov attractor is hyperchaotic and in which it has only one positive Lyapunov exponent and what exactly response for it. For example, the discrete Shilnikov attractors from the nonholonomic model of Chaplygin top \cite{BorKazSat2016} and Celtic stone have only one positive Lyapunov exponent. On the other hand, such attractors are hyperchaotic in three-dimensional H\'enon maps \cite{GonOvsSimoTur2004}, \cite{GonGonPhD16}, the modified oscillator of Anishchenko-Astakhov and the model under consideration. We suppose that in some cases two-dimensional area expansion near the saddle-focus orbit with two-dimensional unstable manifold is compensated by the volume contraction near some other saddle periodic (quasiperiodic) orbits that also belongs to the attractor but which has only one-dimensional unstable manifold, and thus, two-dimensional areas are contracted in their neighborhood. Since Lyapunov exponents are average characteristic of the attractor, only one Lyapunov exponent can become positive in this case.

\section{Multistability in the dynamics of two bubbles} \label{sec:multistability}

\begin{figure}[!ht]
\center{\includegraphics[width=0.4\linewidth]{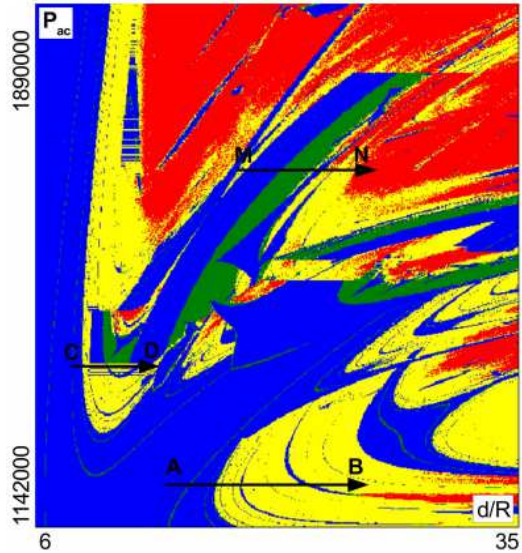} }
\caption{{\footnotesize Another leaf of the chart of Lyapunov exponents for $\omega = 2.87 \cdot 10^7$ s$^{-1}$. Comparsion to Fig. \ref{fig:clock} demonstrates signifacnt presence of multistability.
}}
\label{fig:2nd_leaf}
\end{figure}


Now let us consider several routes, mentioned in Fig. \ref{fig:clock}a, for which multistability is presented, namely, routes AB, CD, MN. First of all, we provide one more two-dimensional map in Fig. \ref{fig:2nd_leaf}, corresponding to the second leaf of the two-dimensional map presented in Fig. \ref{fig:clock}a. Comparison of Fig.-s \ref{fig:clock}a and \ref{fig:2nd_leaf} shows that there a lot of substantially multistable areas in this parameters region (not even taking into consideration the possibility of coexistence of different attractors of the same type, for example two different periodic attractors,  which are not reflected in this kind of maps).

\begin{figure}[!hb]
\center{\includegraphics[width=0.9\linewidth]{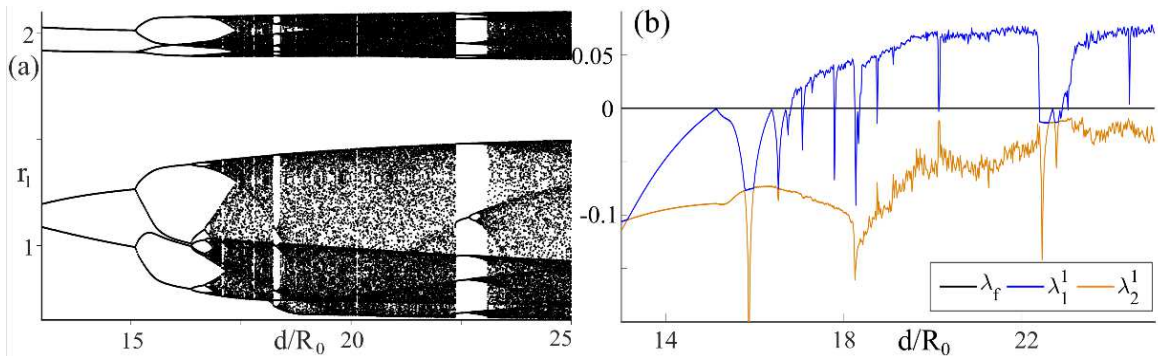} }
\caption{{\footnotesize (a), (b) bifuraction tree and graph of two largest Lyapunov exponents for the bifurcation sequence of the synchronous attractor on the path AB:  $13 \cdot R_0 < d  < 25 \cdot R_0$, $P_{ac} = 1.2$ MPa. Compare to Fig. \ref{fig5a}.
}}
\label{fig6}
\end{figure}


Let us start the one-parameter analysis from the path AB lying in the following parameters interval:$P_{ac} = 1.2$ MPa, $13  R_0 < d < 25  R_0$. We have already shown that one of the attractors existing at the point A (the asynchronous one) goes through the scenario described in the previous section and becomes hyperchaotic, see Fig. \ref{fig5a}. However there exists also a synchronous periodic attractor at the point A, which, on the same path, goes through a Feigenbaum's cascade of the period-doubling bifrucations and becomes synchronous chaotic with one positive Lyapunov exponent. In Fig. \ref{fig6} we present the bifurcation tree and the graph of two largest Lyapunov exponents corresponding to the transition of the synchronous periodic to the synchronous chaotic oscillations on the route AB. This also provides an example of coexisting of a synchronous chaotic oscillations with a hyperchaotic attractor (compare Fig.-s \ref{fig5a} and \ref{fig6}).

\begin{figure}[!ht]
\center{\includegraphics[width=0.9\linewidth]{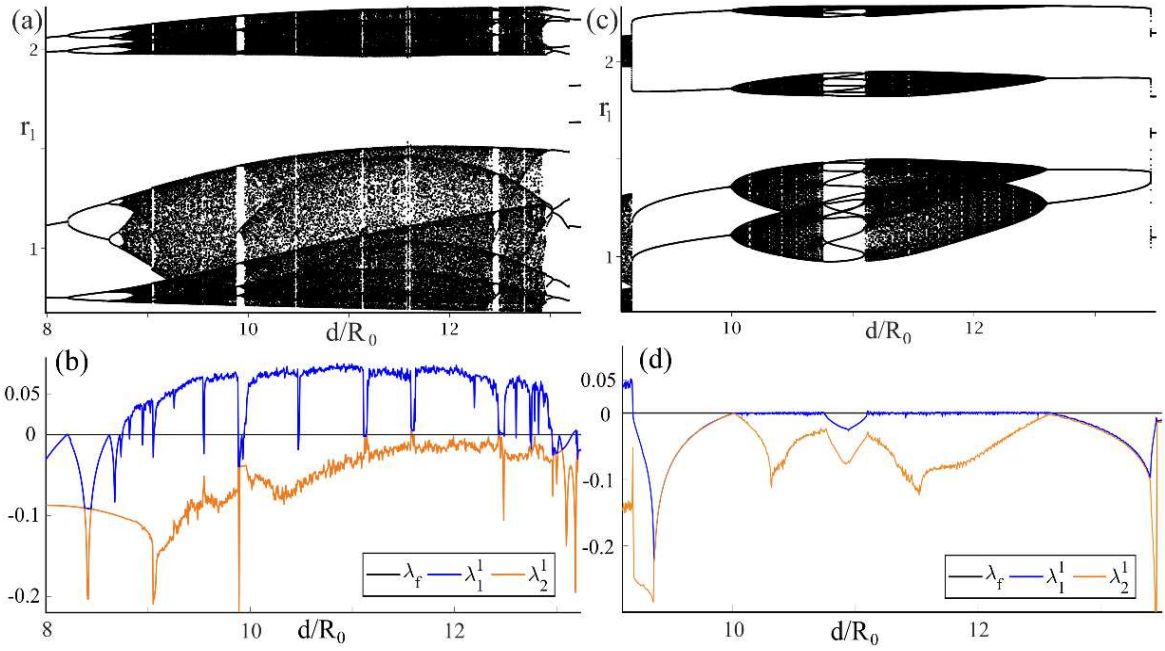} }
\caption{{\footnotesize Route CD:  $8 \cdot R_0 < d  < 13.5 \cdot R_0$, $P_{ac} = 1.4$ MPa. (a), (b) bifurcation tree and two largest Lyapunov exponents, corresponding to the synchronous attractor. (c), (d) bifurcation tree and two largest Lyapunov exponents, corresponding to the asynchronous attractor.
}}
\label{fig7}
\end{figure}


Route CD corresponds to $P_{ac} = 1.4$ MPa and the following $d/R_0$ interval: $8 < d/R_0 < 13.5$, see Fig. \ref{fig7}. No hyperchaotic attractors are presented on this route, but several examples of multistability: coexistence of quasiperiodic oscillations with periodic, quasiperiodic and chaotic, asynchronous periodic and synchronous chaotic, quasiperiodic and synchronous periodic.

\begin{figure}[!ht]
\center{\includegraphics[width=1.0\linewidth]{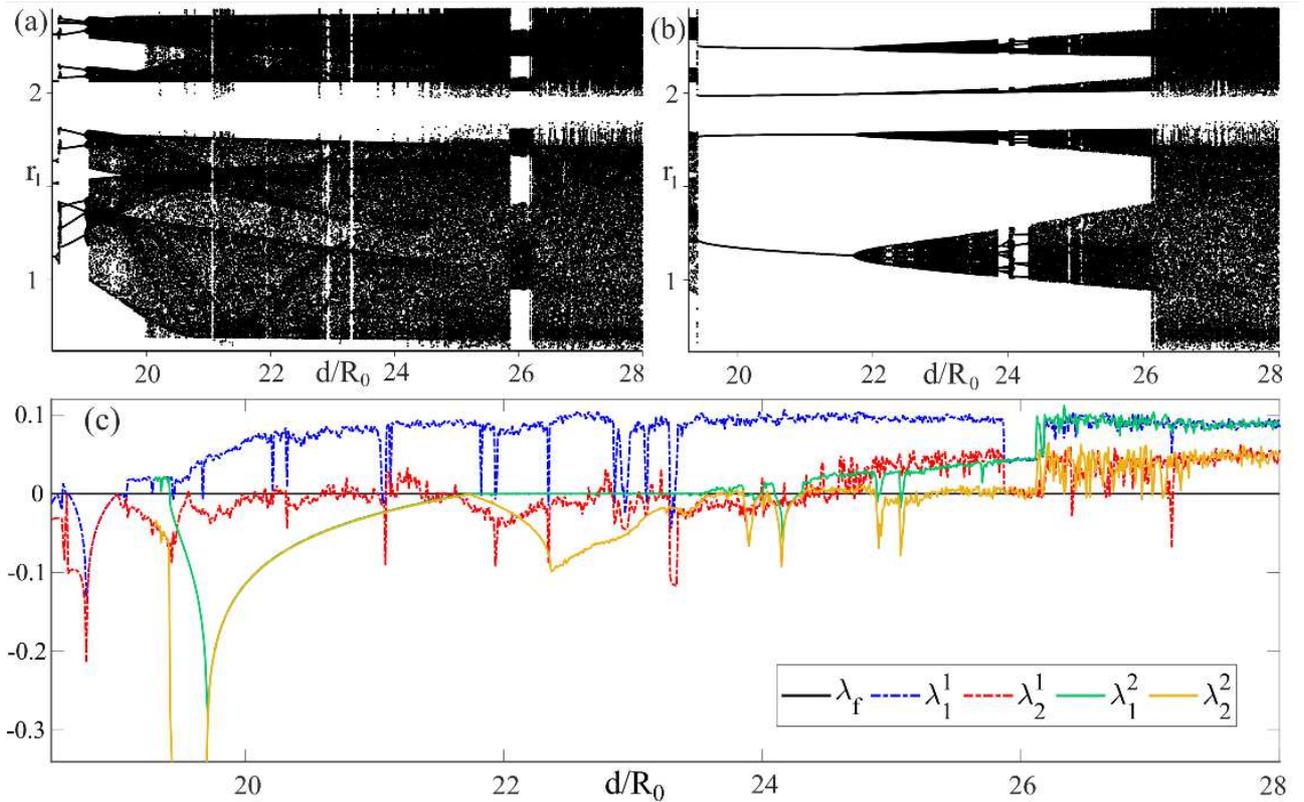} }
\caption{{\footnotesize (a), (b), bifuraction trees corresponding to different attractors and (c) graph of two largest Lyapunov exponents for both the attractors on the route MN:  $18.5 \cdot R_0 < d  < 28 \cdot R_0$, $P_{ac} = 1.68$ MPa. Lyapunov exponents $\lambda^1_1$ and $\lambda^1_2$ correspond to the attractor associated with plot (a), $\lambda^2_1$ and $\lambda^2_2$ -- to the attractor associated with plot (b). A variety of types of multistability can be observed on a single path.
}}
\label{fig8}
\end{figure}


A lot of complicated behavior can be observed on the path MN. For this route we fix $P_{ac} = 1.68$ MPa and vary $d$ in the following interval: $18.5 < d/R_0 < 28$. For this path we provide two bifurcation trees which correspond to different attractors in Fig.-s \ref{fig8}a,b and present the merged picture of the Lyapunov spectra for both attractors in Fig. \ref{fig8}c. We do that in order to explicitly show coexistence of attractors of different types. In Fig. \ref{fig8}c, Lyapunov exponents $\lambda^1_1$ and $\lambda^1_2$ are two largest exponents of the attractor corresponding to Fig. \ref{fig8} a, and $\lambda^2_1$ and $\lambda^2_2$ are the largest exponents of the attractor corresponding to Fig. \ref{fig8}c. $\lambda_f ==0$ is the referent exponent, which is always zero and the same for both attractors. On the right side of the diagram in Fig.-s \ref{fig8}a,c one can observe an abrupt shift from hyperchaos to chaos occuring at $d/R_0 \approx 25.94$. This represents the moment when the chaotic attractor, corresponding to the one in diagram Fig. \ref{fig8}b, remains the only attractor in the system, and the system switches to it. In the remaining interval $25.94 < d/R_0 < 28$, the Lyapunov spectra overlap and the bifurcation diagrams in Fig. \ref{fig8}a,b look identical, because they represent the same attractor. The same feature is observed at $d/R_0 \approx 19.42$. The Lyapunov exponents $\lambda^2_1$ and $\lambda^2_2$ leap at this point and begin to overlap with $\lambda^2_1$, $\lambda^2_2$ to the left of it, see the left side of Fig. \ref{fig8}c. It represents the shift from the disappearing stable limit cycle (the very left side of the bifurcation tree Fig. \ref{fig8}b) to the chaotic attractor corresponding to Fig. \ref{fig8}a. For all the lower values of $d$ these attractor (and their spectra) coincide and we don`t draw all the exponents further left.
 Thus the two attractors discussed here coexist in the interval $19.42 < d/R_0 < 25.94$. A lot of different kinds of multistability can be observed in this interval with these two attractors. In Fig. \ref{fig8}c one can see the following types of coexisting attractors: chaotic with periodic, chaotic with quasiperiodic, hyperchaotic with periodic, hyperchaotic and quasiperiodic, hyperchaotic and chaotic.

Note that the the structure of the border between quasiperiodic and hyperchaotic regimes around the route MN in Fig. \ref{fig:2nd_leaf} looks very similar to the same border around the route AB in Fig. \ref{fig:clock}a. The structure of the graph of the Lyapunov exponents $\lambda^2_1, \lambda^2_2$ in Fig. \ref{fig8}c also looks similar to those in Fig. \ref{fig5a}b and in Fig. \ref{fig5}b. This leads us to the conclusion that it is highely likely that the onset of the hyperchaotic attractor corresponding to Fig. \ref{fig8}b happens by the same scenario that was described in Sec. \ref{sec:hyperchaos}.

\section{Conclusion}

In this work we have studied the nonlinear dynamics of two encapsulated interacting gas bubbles in a liquid. We have showed that the oscillations of bubbles can be periodic, quasyperiodic, chaotic and hyperchaotic. Moreover, we have observed multistability phenomenon in a wide region of the control parameters, which makes bubbles' dynamics even more complicated. We believe that both quasyperiodic and hyperchaotic oscillations along with multistability phenomenon are reported for the first time.


Concernring the onset of chaotic dynamics, we have studied typical roots to chaos and hyperchaos in system \eqref{eq:eq1}. We have demonstrated that simple chaotic attractors (which only one positive Lyapunov exponent) occur either via the Feigenbaum's cascade of period-doubling bifurcations or by the Afraimovich--Shilnikov scenario of the destruction of invariant tori. On the other hand, for the onset of hyperchaotic oscillations we propose a new scenario which is based on the appearance of a discrete Shilnikov attractor containing a saddle-focus periodic orbit with its two-dimensional unstable manifold. Orbits on this attractor can pass arbitrary close to this saddle-focus orbit, where two-dimensional areas are expanded. As a result, two Lyapunov exponents become positive i.e. a hyperchaotic attractor is born. As we know the proposed scenario gives one of few known explanations of the emergency of hyperchaotic behavior. Moreover, we believe that this scenario may be typical for other multidimensional systems demonstrating transition to hyperchaos via the destruction of a torus.

We have also studied multistability phenomenon in system \eqref{eq:eq1}. We have showed that various types of attractors, both synchronous and asynchronous, regular (periodic or quasiperiodic) and chaotic, and even hyperchaotic can co-exist at the same values of the control parameters. In particular, synchronous periodic regimes can coexist with asynchronous periodic, quasiperiodic, asynchronous chaotic or hyperchaotic states, and even coexistence of several sycnhronous regimes is possible (for example, two different synchronous periodic limit cycles). We have demonstrated that in multistable states hyperchaotic regimes can coexist with regular and chaotic (both synchronous and asynchonous) ocillations, as well as with asynchronous qusiperiodic regimes.

As far as applications are concerned, it is known \cite{Hoff,Carroll2013} that chaotic oscillations of bubbles can be beneficial for blood flow visualization. Thus, we believe that the regions of the control parameters where either one chaotic or hyperchaotic attractor exists or both of these attractors co-exist, may be recommended for this type of applications. On the other hand, the regions of the control parameters where different types of attractors exist (e.g. periodic and quasipariodic) should be avoided in applications, since in this case the dynamics of bubbles becomes virtually unpredictable due to the fact that small perturbations in the initial conditions or control parameters may lead to a substantial change in bubbles' acoustic response.

\section*{Acknowledgments}
This paper (except Section 3) was supported by RSF grant 17-71-10241. The results in Section 3 were supported by RFBR grant 18-31-20052. A.K also thanks Basic Research Program at NRU HSE in 2019 for support of scientific researches.

\end{document}